\def \NN {\mathbb N}
\def \CC {\mathbb C}

\def \RR {\mathbb R}

\def \epsilon{\varepsilon}

\def \I  {{\mathcal I}}
\def \J  {{\mathcal J}}
\def \LL {{\mathcal L}}

\def \S  {{\mathcal S}}

\def \X {{\mathcal X}}
\def \d {\text{d}}

\def \ep {\epsilon}
\def \ga {\gamma}

\def \si {\sigma}

\renewcommand{\S}{{\mathcal S}}

\documentclass[12pt,reqno]{amsart}

\usepackage{amsmath,amssymb}

\newtheorem{lem}{Lemma}
\usepackage{url} 

\numberwithin{equation}{section}

\begin{document}
\title[]{Twists by Dirichlet characters and polynomial Euler products of $L$-functions, II}
\author[]{J.Kaczorowski \lowercase{and} A.Perelli}
\maketitle

\hfill {\it In memory of Eduard Wirsing}

\medskip
{\bf Abstract.} In a previous paper we proved that if an $L$-function  $F$ from the Selberg class has degree $2$, its conductor $q_F$ is a prime number and  $F$ is weakly twist-regular at all primes $p\neq q_F$, then $F$ has a polynomial Euler product. In this paper we extend this result to $L$-functions of degree 2 with square-free conductor $q_F$, which are weakly twist-regular at all primes $p\nmid q_F$.



\section{Introduction} 

\smallskip
In \cite{Ka-Pe/twistI} we proved, among other things, that if an $L$-function  $F$ from the Selberg class $\S$ has degree $2$, its conductor $q_F$ is a prime number, and  $F$ is weakly twist-regular at all primes $p\neq q_F$, then $F$ has a polynomial Euler product. In this paper we keep the notation from \cite{Ka-Pe/twistI} but, for the reader's convenience, we recall some basic definitions in Section \ref{sec-def}. The aim of this paper is to extend the above result as follows.

\medskip
{\bf Theorem 1.} {\it Let $F\in \S$ be of degree $2$ and its conductor $q_F$ be square-free. If $F$ is weakly twist-regular at all primes $p \nmid q_F$, then $F$ has a polynomial Euler product.}

\medskip
This is not a straightforward generalization; indeed, apart from the use of Theorem 2 below, a non-trivial extension of the method in \cite{Ka-Pe/twistI} is necessary. On the other hand, the present method cannot settle the problem in full generality, i.e. for all integer conductors, and some new ideas will probably be needed to prove the full result.

\smallskip
One of the main tools in the proof is the following transformation formula for linear twists of $L$-functions from the extended Selberg class, which is of independent interest.

\medskip
{\bf Theorem 2.} {\it Let $F\in \S^\sharp$ with $d_F=2$, and let $\alpha>0$. Then for every integer $K>0$ there exist polynomials $Q_0(s),...,Q_K(s)$, with $Q_0(s)\equiv 1$, such that
\begin{equation}
\label{eq:Falpha} 
F(s,\alpha) = -i\omega_F^* (\sqrt{q_F}\alpha)^{2s-1+2i\theta_F} \sum_{\nu=0}^K \big(\frac{iq_F\alpha}{2\pi}\big)^\nu Q_\nu(s) \overline{F}\big(s+\nu+2i\theta_F,-\frac{1}{q_F\alpha}\big) + H_K(s,\alpha).
\end{equation}
Here $H_K(s,\alpha)$ is holomorphic for $-K+\frac12<\sigma<2$ and satisfies
\begin{equation}
\label{HK}
H_K(s,\alpha) \ll (|s|+1)^{2K + A}
\end{equation}
with a certain constant $A=A(F,\alpha)>0$. Moreover, $\deg Q_\nu=2\nu$ and
\begin{equation}
\label{eq:Q1}
Q_\nu(s) \ll \frac{(A(|s|+1))^{2\nu}}{\nu!} \hskip1.5cm \text{for} \ 1\leq \nu\leq \min(|s|,K).
\end{equation}
}

\medskip
This should be compared with Theorem 1.2 in \cite{Ka-Pe/2015}. Apart from the value $2\theta_F$ in place of $\theta_F$ in \eqref{eq:Falpha}, due to a slight change in the definition of $\theta_F$ compared with \cite{Ka-Pe/2015} (see next section), the main difference is in the estimate for the size of $H_K(s,\alpha)$ and in the range for $s$ in which it holds. Precisely, Theorem 1.2 in  \cite{Ka-Pe/2015} states that 
\begin{equation}
\label{HK1}
H_K(s,\alpha)\ll (AK)^K \hskip1.5cm \text {for $-K+\frac12<\sigma<2$, $|s|\leq 2K$}. 
\end{equation}
In the proof of Theorem 1, we consider shifts of $L$-functions of the form $F(s+i\tau)$ with $\tau\to\infty$ and hence reasonable control upon the size of $H_K(s,\alpha)$ is needed in half-planes rather than in discs. If fact, for $\tau\to\infty$, (\ref{HK1}) gives an estimate which is far too weak for the proof of Theorem 1. In contrast, (\ref{HK}) secures a polynomial growth which is exactly what is needed. In principle, the main structure of the proof of Theorem 2 is the same as that of Theorem 1.2 in \cite{Ka-Pe/2015}. Nevertheless, the new situation where $s$ is in a half-plane forces the introduction of significant technical changes in the proof.

\medskip
{\bf Acknowledgements}.  We wish to warmly thank the referee for carefully reading our manuscript and for pointing out many inaccuracies. This research was partially supported by the Istituto Nazionale di Alta Matematica, by the MIUR grant PRIN-2017 {\sl ``Geometric, algebraic and analytic methods in arithmetic''} and by grant 2021/41/BST1/00241 {\sl ``Analytic methods in number theory''}  from the National Science Centre, Poland.

\medskip
\section{Definitions}\label{sec-def}

\smallskip
Throughout the paper we write $s=\si+it$, $e(x) = e^{2\pi ix}$ and $\overline{f}(s)$ for $\overline{f(\overline{s})}$. The extended Selberg class $\S^\sharp$ consists of non-identically vanishing Dirichlet series
\[
F(s) = \sum_{n=1}^\infty \frac{a(n)}{n^s}
\] 
absolutely convergent for $\si>1$, such that $(s-1)^mF(s)$ is entire of finite order for some integer $m\geq0$, and satisfying a functional equation of type
\[
F(s) \gamma(s) = \omega \overline{\gamma}(1-s) \overline{F}(1-s),
\]
where $|\omega|=1$ and the $\gamma$-factor
\[
\gamma(s) = Q^s\prod_{j=1}^r\Gamma(\lambda_js+\mu_j) 
\]
has $Q>0$, $r\geq0$, $\lambda_j>0$ and $\Re(\mu_j)\geq0$. Note that the conjugate function $\overline{F}$ has conjugated coefficients $\overline{a(n)}$. The Selberg class $\S$ is the subclass of $\S^\sharp$ of the functions satisfying the Ramanujan conjecture $a(n) \ll n^\ep$ and with an Euler product of the form
\[
F(s)=\prod_p F_p(s), \quad \text{where} \quad F_p(s) = \sum_{k=0}^{\infty} \frac{a(p^k)}{p^{ks}}
\]
satisfies
\begin{equation}
\label{logeffepi}
\log F_p(s) = \sum_{k=1}^{\infty}\frac{b(p^k)}{p^{ks}} \quad \text{with $ b(p^k)\ll p^{\vartheta k}$ for a certain $\vartheta <1\slash 2$.}
\end{equation}
Note that the series in \eqref{logeffepi} is absolutely convergent for $\si>\vartheta$ and hence 
\begin{equation}
\label{effepi}
\text{$F_p(s)$ is holomorphic and bounded away from 0 for $\si>\vartheta'$ for some $\vartheta'<1/2$.}
\end{equation}
We say that $F\in \S$ has a polynomial Euler product if for every prime $p$ there exist $\partial_p\in\NN$ and $\alpha_{j,p}\in\CC$ such that
\[ 
F_p(s) = \prod_{j=1}^{\partial_p} \left( 1-\frac{\alpha_{j,p}}{p^s}\right)^{-1}.
\]
In such a case, it follows from the Ramanujan conjecture that $|\alpha_{j,p}|\leq 1$; see p.448-449 of \cite{Ka-Pe/2015}. 

\smallskip
The twist by a Dirichlet character $\chi$ (mod $q$) is defined for $\si>1$ as
\[ 
F^{\chi}(s) = \sum_{n=1}^{\infty}\frac{a(n)\chi(n)}{n^s}.
\] 
We say that $F\in \S^\sharp$ is weakly twist-regular at $p$ if for every primitive Dirichlet character $\chi$ (mod $p^f$) with $1\leq f\leq m_{q_F}(p)$, where $m_{q_F}(p)$ is the order of $p$ (mod $q_F$), the twist $F^\chi$ belongs to $\S^\sharp$ and has the same degree as $F$. Moreover, the linear twist of $F\in\S^\sharp$ is defined for $\si>1$ as
\[
F(s,\alpha) = \sum_{n=1}^\infty \frac{a(n)}{n^s} e(-n\alpha)
\]
with $\alpha\in\RR$.

\smallskip
Degree $d_F$, conductor $q_F$, root number $\omega_F$ and $\xi$-invariant $\xi_F$ of $F\in\S^\sharp$ are defined by
\[
\begin{split}
d_F =2\sum_{j=1}^r\lambda_j&, \qquad q_F= (2\pi)^dQ^2\prod_{j=1}^r\lambda_j^{2\lambda_j}, \\
\omega_F=\omega \prod_{j=1}^r \lambda_j^{-2i\Im(\mu_j)}&,  \qquad \xi_F = 2\sum_{j=1}^r(\mu_j-1/2)= \eta_F+ id_F\theta_F
\end{split}
\]
with $\eta_F,\theta_F\in\RR$. We also write
\[
\omega_F^* = \omega_F e^{-i\frac{\pi}{2}(\eta_F+1)}\big(\frac{q_F}{(2\pi)^{2}}\big)^{i\theta_F} \quad \text{and} \quad \tau_F=\max_{j=1,\dots,r}\big|\frac{\Im(\mu_j)}{\lambda_j}\big|,
\]
while $m_F$ denotes the order of the pole of $F$ at $s=1$. In $\omega_F^*$, and in other definitions below, we changed $\theta_F$ to $2\theta_F$ compared to the corresponding definitions in \cite{Ka-Pe/2015}. This is due to the above slightly different definition of the $\xi$-invariant and the fact that we are considering functions of degree $d_F=2$. The $H$-invariants of $F$ are defined as
\[
H_F(n) = 2\sum_{j=1}^r\frac{B_n(\mu_j)}{\lambda_j^{n-1}} \hskip1.5cm n=0,1,...
\]
where $B_n(x)$ is the $n$-th Bernoulli polynomial. Note that $H_F(0) = d_F$ is the degree and 
$H_F(1)$  is the $\xi$-invariant. We refer to our survey papers \cite{Kac/2006},\cite{Ka-Pe/1999b},\cite{Per/2005},\cite{Per/2004} for further definitions, examples and the basic theory of the Selberg class. 

\smallskip
As in \cite{Ka-Pe/2015}, for $\nu,\mu=1,2,...$ we define the polynomials $R_\nu(s) = R_{\nu,F}(s)$ and $V_\mu(s)=V_{\mu,F}(s)$ as
\begin{equation}
\label{3-0bis}
\begin{split}
R_\nu(s) &= B_{\nu+1}(1-2s-2i\theta_F) + B_{\nu+1}(1) \\
&+ \frac12  \sum_{k=0}^{\nu+1}{\nu+1\choose k}\big((-1)^\nu H_F(k)s^{\nu+1-k} - \overline{H_F(k)}(1-s)^{\nu+1-k}\big)
\end{split}
\end{equation}
and
\begin{equation}
\label{3-1}
V_\mu(s) = (-1)^\mu\sum_{m=1}^\mu \frac{1}{m!} \sum_{\substack{\nu_1\geq 1,...,\nu_m\geq 1\\ \nu_1+...+\nu_m=\mu}} \prod_{j=1}^m\frac{R_{\nu_j}(s)}{\nu_j(\nu_j+1)},
\end{equation}
respectively. We also define $Q_0(s)\equiv1$ and, for $\nu=1,2,...$, the polynomials $Q_\nu(s)$ appearing in Theorem 2 by means of the formula
\[
\exp\big(\sum_{\nu=1}^\infty \frac{(-1)^\nu R_\nu(s)}{\nu(\nu+1)} \frac{1}{(w+2s-1+2i\theta_F)^\nu}\big) \approx 1 + \sum_{\nu=1}^\infty \frac{Q_\nu(s)}{(w-1)\cdots(w-\nu)},
\]
and the coefficients $C_{\mu,\ell}$, $\ell\geq\mu\geq1$, by
\begin{equation}
\label{ref-1}
\frac{1}{w^\mu} \approx \sum_{\ell=\mu}^\infty \frac{C_{\mu,\ell}}{(w-1)\cdots(w-\ell)}.
\end{equation}
Moreover, we define the coefficients 
$A_{\mu,\nu}(s)$ ($\nu\geq \mu\geq 1$) by
\begin{equation}
\label{3-3}
\frac{1}{(w+2s-1+2i\theta_F)^\mu} \approx  \sum_{\nu=\mu}^\infty \frac{A_{\mu,\nu}(s)}{(w-1)\cdots(w-\nu)}.
\end{equation}
Here $\approx$ means asymptotic expansion as $w\to\infty$, i.e. cutting the sum on the left hand side at $\nu=N$ introduces an error of size $O_s(1/|w|^{N+1})$. We refer to Lemmas \ref{lemwmu} and  \ref{lemwmu2} below for a more precise meaning of  \eqref{ref-1} and \eqref{3-3}.

\smallskip
We write $w=u+iv$ and, for a given $s$, define the contour $\LL(s)$ as follows:
\[
\LL(s) = \LL_{-\infty}(s) \cup \LL_{\infty}(s)
\]
where
\[
\LL_{-\infty}(s) = (-\sigma+c_0-i\infty,-\sigma+c_0+it_0] \cup [-\sigma+c_0+it_0,-\sigma-c_0+it_0]
\]
\[
\LL_\infty(s) = [-\sigma-c_0+it_0,-\sigma-c_0+i\infty).
\]
Here $t_0=t_0(s) = c_1(|s|+1)^2$ and $c_0,c_1>0$ are sufficiently large constants depending on $F$ to be chosen later on. Observe a significant difference in the present choice of $t_0$ compared to the analogous choice in \cite{Ka-Pe/2015}; here $t_0$ is much larger, and this is important in the proof of Theorem 2. Moreover, we denote by $\LL^*_{-\infty}(s)$ the half-line $1-2s-2i\theta_F - \LL_\infty(s)$ taken with the positive orientation, hence
\[
\LL_{-\infty}^*(s) = (1-\sigma+c_0-i\infty, 1-\sigma+c_0-it_0^*]
\]
with $t_0^*=t_0^*(s) =t_0 +2t+2\theta_F$. Further, we let
\[
\LL_\infty^*(s)=[1-\sigma+c_0-it_0^*, N+1] \cup [N+1,N+1+i\infty),
\]
where the positive integer $N$ will be chosen later on (see \eqref{P6bis} below), and write
\[
\LL^*(s) = \LL^*_{-\infty}(s) \cup \LL_\infty^*(s).
\]

\smallskip
We shall also use the notation
\[
G(s,w) = \frac{(2\pi)^{1-r}}{\Gamma(1-w)} \prod_{j=1}^r \Gamma(\lambda_j(1-s-w)+\bar{\mu}_j) \Gamma(1-\lambda_j(s+w)-\mu_j)
\]
and
\[
S(s,w) = \frac{2^{r-1}}{\sin\pi w} \prod_{j=1}^r \sin\big(\pi(\lambda_j(s+w)+\mu_j\big).
\]

\smallskip
Finally, $A, B, c,c',...$ will denote positive constants, possibly depending on $F$ (also via a dependence on the above constants $c_0,c_1$), not necessarily the same at each occurrence. The constants in the $\ll$- and $O$-symbols may also depend on $F(s)$ (again, also via $c_0,c_1$).

\medskip
\section{Lemmas}

\begin{lem}\label{lem:sum}
Let $F\in\S$ and $q_F$ be square-free. Then for $\si>1$ we have
\begin{equation}\label{eq:sum}
\sum_{\substack{1\leq a\leq q_F\\ (a,q_F)=1}} \sum_{n=1}^{\infty} \frac{a(n)}{n^s} e(-\frac{a}{q_F}n)
= F(s) \sum_{d|q_F} d \mu(\frac{q_F}{d}) \prod_{p|d}\left(1-F_p(s)^{-1}\right).
\end{equation}
\end{lem}

{\it Proof.} For simplicity we write $q$ in place of $q_F$. Recall the well-known Kluyver's formula for the Ramanujan sum 
\[ 
c_q(n):= \sum_{\substack{1\leq a\leq q\\ (a,q)=1}} e(an\slash q)= \sum_{d|(q,n)}
\mu(\frac{q}{d}) d.
\]
Since $c_q(n)\in\RR$, the left hand side of (\ref{eq:sum}) equals
\begin{equation}\label{eq:sum1}
\sum_{n=1}^{\infty}\frac{a(n)}{n^s} c_q(n) = \sum_{d|q}\mu(\frac{q}{d})d^{1-s}
\sum_{n=1}^{\infty}\frac{a(dn)}{n^s}.
\end{equation}
Now we compute the inner sum for a generic  $d|q$, $d=p_1\ldots p_k$. We have
\begin{eqnarray*}
\sum_{n=1}^{\infty}\frac{a(dn)}{n^s}
 &=& 
\sum_{\nu_1=0}^{\infty}\ldots\sum_{\nu_k=0}^{\infty}
\sum_{\substack{n\geq 1\\ p_j^{\nu_j}||n \\
(1\leq j\leq k)}}\frac{a(dn)}{n^s}\\
&=&\sum_{\nu_1=0}^{\infty}\ldots \sum_{\nu_k=0}^{\infty}
\left(\prod_{j=1}^k \frac{a(p_j^{\nu_j+1})}{p_j^{\nu_js}}\right)
\sum_{(n,p_1\ldots p_k)=1}\frac{a(n)}{n^s}\\
&=& \sum_{\nu_1=0}^{\infty}\ldots \sum_{\nu_k=0}^{\infty}
\left(\prod_{j=1}^k \frac{a(p_j^{\nu_j+1})}{p_j^{\nu_js}}\right)
\left(\prod_{j=1}^k F_{p_j}(s)^{-1}\right) F(s)\\
&=&F(s)d^s\prod_{j=1}^k \left(1-F_{p_j}(s)^{-1}\right).
\end{eqnarray*}
Inserting this into (\ref{eq:sum1}) we obtain (\ref{eq:sum}), and the proof is complete. \qed

\medskip
For $F$ in the Selberg class $\S$ we denote by $N_F(\si,T)$ the number of non-trivial zeros $\beta+i\gamma$ of $F$ in the rectangle $\beta> \si$, $|\ga|\leq T$.

\begin{lem}\label{lem:density}   Let $F\in \S$ with $d=2$. Then for every $\ep>0$ and any fixed $\si>1/2$ we have
\[
N_F(\si,T) \ll T^{3/2-\si+\ep}.
\]
\end{lem}

{\it Proof.} See p.474-475 of \cite{Ka-Pe/2015}. \qed

\begin{lem}\label{lem:estimate} Let $F\in \S^{\sharp}$ with $d=2$. Then there exists a positive constant $T_0=T_0(F)$ such that
\[
F(s) \asymp \left(\frac{q_F}{(2\pi e)^2}\right)^{|\si|} |s|^{2|\si|+1}
\]
uniformly for $|t|\geq T_0$ and $\si\leq-1$. 
\end{lem}

{\it Proof.} This is a  refined version of Lemma 2.1 in \cite{Ka-Pe/2015}. The proof follows from the functional equation of $F$ and the Stirling formula; see also Lemma 3 in \cite{Ka-Pe/twistI}. \qed

\begin{lem}\label{lem:cchi}
Let  $F\in \S$, $p$ be a prime number and $\si>1$. Then there exist coefficients $c(\chi,p)$, where $\chi$ runs over the Dirichlet characters $\chi$ {\rm (mod $p$)}, such that
\[
F(s,1/p) = \sum_{\substack{\chi(\bmod p) \\ \chi\neq \chi_0}} c(\chi,p) F^{\chi}(s) + \left(1- \frac{p}{p-1}F_p(s)^{-1}\right) F(s).
\]
\end{lem}

{\it Proof.} See equation (2.6) of \cite{Ka-Pe/2015}, observing that $F(s,1/p) = F(s,-(p-1)/p)$. \qed

\begin{lem}\label{lem:tauk}
Let $F\in \S^{\sharp}$, $q$ be a positive integer and let $\theta<1\slash 2$ be fixed. Moreover, for every prime $p|q$, let $\varepsilon_p$ be a complex number with $|\ep_p|=1$.  Then there exist two positive constants  $a$ and $b$ (depending only on $\theta$) and a sequence of positive numbers $\tau_k\to\infty$, $k\geq1$, satisfying the following two conditions

i) for every prime $p|q$ we have $|p^{-i\tau_k} - \varepsilon_p|<1\slash k$,

ii) $|F(\sigma +i t)|\geq \tau_k^{-b}$  uniformly for $-1\leq\si\leq \theta$ and $|t-\tau_k|\leq \tau_k^a.$ 

\end{lem}

{\it Proof.} The numbers $\log p$, $p|q$ are linearly independent over $\mathbb Q$. Thus by a well known version of the classical Kronecker approximation theorem, for every $k\geq 1$ there exists a relatively dense set of numbers $\tau$ such that 
\[
|p^{-i\tau}-\varepsilon_p|< \frac{1}{k}  \quad \text{for all} \ p|q.
\]
Thus using Lemma \ref{lem:density} we see that there exists a solution $\tau_k\geq 2\max(k, T_0(F))$, where $T_0(F)$ is the constant appearing in Lemma \ref{lem:estimate}, of this system of inequalities such that $F(\si+it)\neq 0$ for 
\[
-1\leq  \si<(1+2\theta)\slash 4\ \quad \text{ and} \quad   |t-\tau_k|\leq \tau_k^a,
 \]
where $a$ is any fixed positive number less than $(1-2\theta)\slash 4$. For $s$ in such a range we use the well known formula 
\[
-\frac{F'}{F}(s) = \sum_{|\gamma - t|<1}\frac{1}{s-\rho} + O(\log \tau_k), \quad \text{therefore} \quad  
-\frac{F'}{F}(s) = O(\log\tau_k).\]
Hence, recalling Lemma \ref{lem:estimate},   we have
\[
\log|F(\si+it)| = \log|F(-1+it)| +\Re \int_{-1}^{\si} \frac{F'}{F}(u+it)\, du\geq - b\log \tau_k
\]
for $-1\leq\si\leq \theta$ and $ |t-\tau_k|\leq \tau_k^a$,  and the lemma follows. \qed

\begin{lem}\label{lemR}  Let $R_{\nu}(s)$ be as in \eqref{3-0bis} and let $c\geq 1$. For $s\in\CC$ and $1\leq \nu\leq c(|s|+1)$ we have
\[
R_\nu(s) \ll \big(c'(|s|+1)\big)^{\nu+1}
\]
with a suitable $c'\geq c+1$.
\end{lem}

{\it Proof.} This is Lemma 3.8 in \cite{Ka-Pe/2015}, with an explicit lower bound for $c'$. \qed

\medskip
Let
\begin{equation}
\label{eqnew1}
V_{\mu,N}(s) = (-1)^\mu \sum_{m=1}^\mu \frac{1}{m!}  \sum_{\substack{1\leq \nu_1\leq N,...,1\leq \nu_m\leq N\\ \nu_1+...+\nu_m=\mu}} \prod_{j=1}^m\frac{R_{\nu_j}(s)}{\nu_j(\nu_j+1)}.
\end{equation}
Recalling definition \eqref{3-1} we see that $V_{\mu,N}(s) = V_{\mu}(s)$ for $\mu\leq N$.

\begin{lem} \label{lemVmuN}  
For $s\in\CC$ and $1\leq \mu\leq c(|s|+1)$ we have
\[
V_{\mu,N}(s)\ll \frac{(c'(|s|+1)^{2\mu}}{(\mu-1)!},
\]
where $c$ and $c'$ are as in Lemma $\ref{lemR}$. Moreover, for all $\mu\geq 1$
\[
V_{\mu,N}(s)\ll (c'(|s|+1)^{2\mu}.
\]
\end{lem}

{\it Proof.} Thanks to Lemma \ref{lemR}, for every $\mu\geq 1$ we have
\begin{equation}
\label{L3-15/2}
\begin{split}
|V_{\mu,N}(s)| &\leq \sum_{m=1}^\mu \frac{1}{m!}  \sum_{ \nu_1+...+\nu_m=\mu} \prod_{j=1}^m\frac{(c'(|s|+1))^{\nu_j+1}}{\nu_j(\nu_j+1)} \\
&\ll (c'(|s|+1))^\mu \sum_{m=1}^\mu \frac{(c'(|s|+1))^m}{m!}\ll (c'(|s|+1))^{2\mu}.
\end{split}
\end{equation}
 Moreover, recalling that $c'\geq c+1$, for $1\leq \mu\leq c(|s|+1)$ we have
\[ 
 \sum_{m=1}^\mu \frac{(c'(|s|+1))^m}{m!}\leq \mu \max_{1\leq m\leq \mu} \frac{(c'(|s|+1))^m}{m!}=
 \frac{(c'(|s|+1)^{\mu}}{(\mu-1)!},
 \]
 and the lemma follows from this and \eqref{L3-15/2}. \qed
 
\begin{lem} \label{lemwmu} For $w\in\CC$ and integers $\mu,M$ with $1\leq \mu\leq M\leq |w|/2$ we have
\begin{equation}\label{*}
\frac{1}{w^\mu} = \sum_{\ell =\mu}^M \frac{C_{\mu,\ell}}{(w-1)\cdots(w-\ell)} + R_{\mu,M}(w),
\end{equation}
where the coefficients $C_{\mu,\ell}$ are defined by \eqref{ref-1} and
\[
R_{\mu,M}(w) \ll \frac{2^MM!}{(\mu-1)!} \frac{1}{|w(w-1)\cdots(w-M)|}.
\]
\end{lem}

{\it Proof.} This is Lemma 3.13 in \cite{Ka-Pe/2015}. \qed

\begin{lem} \label{lemwmu2} Let $s\in\CC$ and $1\leq \mu\leq\nu\leq N$ be integers and let
\[
A_{\mu,\nu}(s) = \sum_{k=0}^{\nu-\mu} {-\mu\choose k} C_{\mu+k,\nu}(2s-1+2i\theta_F)^k.
\]
Then for $w\in\LL_{-\infty}^*(s)$ and $1\leq N \leq |\sigma| +c$ we have
\[
\frac{1}{(w+2s-1+2i\theta_F)^\mu} = \sum_{\nu=\mu}^N \frac{A_{\mu,\nu}(s)}{(w-1)\cdots(w-\nu)} + 
O\big(\frac{4^NN!|2s-1+2i\theta_F|^{N-\mu+1}}{|w(w-1)\cdots(w-N)|}\big),
\]
where the constant in the $O$-symbol may depend on $c$.
\end{lem}

{\sl Proof.} We follow the proof of Lemma 3.14 in \cite{Ka-Pe/2015}, assuming that the constant $c_1$ in the definition of $\LL_{-\infty}^*(s)$ in Section 2 is sufficiently large and estimating in a different way the error term. Given an integer $P\geq 0$, for $w\in\LL_{-\infty}^*(s)$ we have that
\begin{equation}
\label{L3-14/2}
\frac{1}{(w+\eta)^\mu} = \sum_{k=0}^P {-\mu\choose k} \frac{\eta^k}{w^{\mu+k}} + O\left( \frac{2^{P+\mu}|\eta|^{P+1}}{|w|^{P+\mu+1}}\right),
\end{equation}
where  $\eta = 2s-1+2i\theta_F$. 
 Choose $P=N-\mu$. Since $|w|\geq 2N$, for $w\in\LL_{-\infty}^*$ the error term in \eqref{L3-14/2} is
\begin{equation}
\label{L3-14/3}
\ll \frac{4^N |\eta|^{N-\mu+1}}{|w(w-1)\cdots(w-N)|}.
\end{equation}
Now we apply Lemma \ref{lemwmu}  with $M=N$ and replace $1/w^{\mu+k}$ in (\ref{L3-14/2}) by the main term of (\ref{*}). This produces an error of size 
\begin{equation}
\label{L3-14/4}
\ll \frac{2^NN!}{|w(w-1)\cdots(w-N)|}\sum_{k=0}^{N-\mu}\left|{-\mu\choose k}\right| \frac{|\eta|^k}{(k+\mu)!}  
\ll  \frac{4^NN! |\eta|^{N-\mu}}{|w(w-1)\cdots(w-N)|},
\end{equation}
hence the expansion of $1/(w+\eta)^\mu$ follows from \eqref{L3-14/2}--\eqref{L3-14/4}, and the lemma is proved. \qed

\begin{lem} \label{3.18}   Let $s\in\CC$ and $c\geq 1$. For $1\leq \nu\leq |s|+c$ we have
\[
Q_\nu(s) \ll \frac{(c'(|s|+1))^{2\nu}}{\nu!}
\]
with a suitable $c'=c'(c)>0$.
\end{lem}

{\it Proof.} This is Lemma 3.18 in  \cite{Ka-Pe/2015}. \qed

\medskip
\section{Proof of Theorem 1} 

\smallskip
From Theorem 2 in \cite{Ka-Pe/twistI} we know that $F_p(s)^{-1}$ is a polynomial in $p^{-s}$ for all primes $p\nmid q_F$. The main difficulty lies in proving that the same holds for the primes $p|q_F$.

\smallskip
Let ${\mathbb T}=\{z\in {\mathbb C}: |z|=1\}$ denote the unit circle, $r$ be the number of prime factors of $q_F$$, {\boldsymbol \varepsilon}\in {\mathbb T}^r$ be given and let $\tau_k$ be the sequence appearing in Lemna \ref{lem:tauk}, with the choice $q=q_F$. Let $l_k:=\log \tau_k$ and assume without loss of generality that the $\tau_k$'s are sufficiently large.  We set
\[ 
\Omega_k:=\{s=\sigma+it: \sigma\leq 1/2, |t|\leq l_k^2\},
\]
\[
\Omega_{k,1}:=\{s\in\Omega_k: \sigma\leq -l_k\} \quad \text{and} \quad \Omega_{k,2}:=\{s\in\Omega_k: \sigma\geq -l_k\}.
\]
Moreover, for $\si>1$ let
\begin{equation}\label{eq:Gkdef}
G_k(s):=\frac{1}{F(s+i\tau_k)}\sum_{\substack{ 1\leq a\leq q_F\\ (a,q_F)=1}} F(s+i\tau_k, a/q_F).
\end{equation}
By Lemma \ref{lem:sum} we have
\begin{equation}\label{eg:Gk}
G_k(s)= \sum_{d|q_F}d\mu(\frac{q_F}{d})\prod_{p|d}\left(1- F_p(s+i\tau_k)^{-1}\right)
\end{equation}
hence it follows from \eqref{effepi} that
\begin{equation}\label{eq:Gk1}
G_k(s)\ll 1 \quad \text{for $\sigma>\vartheta'$ with a certain $\vartheta'<1/2$}.
 \end{equation}

\smallskip
As in the proof of Theorem 3 in \cite{Ka-Pe/twistI}, for every $a(\bmod q_F)$, $(a,q_F)=1$, we fix a prime $p_a\equiv a(\bmod q_F)$. Obviously $F(s,a/q_F)=F(s,p_a/q_F)$, thus we apply Theorem 2 with $s\in\Omega_k$, $\alpha=p_a/q_F$ and $K=[|\sigma|]+2$ to obtain that
\begin{eqnarray}\label{eq:tr}
&F&\!\!\!\!\!(s+i\tau_k,a/q_F) = -i\omega_F^*\left(\frac{p_a}{\sqrt{q_F}}\right)^{2s-1+2i(\tau_k+\theta_F)}\nonumber\\ &\times& \sum_{\nu=0}^K\left(\frac{ip_a}{2\pi}\right)^{\nu} Q_{\nu}(s+i\tau_k) \overline{F}(s+\nu + i(\tau_k+2\theta_F), -1/p_a) 
+ H_K(s+i\tau_k,p_a/q_F).
\end{eqnarray}
Moreover, using Lemma \ref{lem:cchi} we have
\begin{equation}\label{eq:tr1}
\overline{F}(s,-1/p_a) = \sum_{\substack{\chi(\bmod p_a)\\\chi\neq\chi_0}} \overline{c({\chi},p_a)} \overline{F^{\chi}}(s) -
\left(1-\frac{p_a}{p_a-1}\overline{F}_{p_a}(s)^{-1}\right) \overline{F}(s).
\end{equation}
From Theorem 2 of \cite{Ka-Pe/twistI} we know that $\overline{F}_{p_a}(s)^{-1}$   is a polynomial in $p_a^{-s}$, and in particular is entire. Since for $s\in \Omega_k$ we have that $\Im(s+\nu + i(\tau_k+2\theta_F))\neq0$, and the possible pole of $F(s)$ and $F^\chi(s)$ is at $s=1$, from (\ref{eq:tr1}) we deduce that $\overline{F}(s+\nu + i(\tau_k+2\theta_F), -1/p_a)$ is holomorphic for $s\in\Omega_k$. Thus from (\ref{eq:Gkdef}),(\ref{eq:tr}) and ii) of Lemma \ref{lem:tauk} we conclude that $G_k(s)$ is holomorphic for $s\in\Omega_k$. 

\smallskip
Now we estimate $G_k(s)$ in this region. From \eqref{eq:Gkdef},\eqref{eq:tr} and \eqref{eq:tr1} we have
\begin{equation}
\label{eq:0}
\begin{split}
G_k(&s) \ll \frac{B^{|\sigma|}}{|F(s+i\tau_k)|} \sum_{0\leq\nu\leq K}|Q_{\nu}(s+i\tau_k)| \\
&\times \max_{\substack{ 1\leq a\leq q_F\\ (a,q_F)=1}} \left(\sum_{\substack{\chi(\bmod p_a)\\\chi\neq\chi_0}}|\overline{F^{\chi}}(s+\nu+i(\tau_k+2\theta_F))| + B^{|\sigma+\nu|}|\overline{F}(s+\nu+i(\tau_k+2\theta_F))|\right) \\ 
&+ \max_{\substack{ 1\leq a\leq q_F\\ (a,q_F)=1}} \frac{|H_K(s+i\tau_k,p_a/q_F)|}{|F(s+i\tau_k)|},
\end{split}
\end{equation}
where here and later on $B>1$ denotes a certain constant, not necessarily the same at each occurrence. Note that we used the fact that $F_{p_a}(s)^{-1}$ is a polynomial in $p_a^{-s}$ and thus it is $\ll B^{|\sigma|}$. Recalling (\ref{eq:Q1}), Lemma \ref{lem:estimate}  and our assumption that $F^\chi$ is a degree 2 function in $\S^\sharp$, we see that for $\si\leq -1$ the contribution of the terms with $\nu\leq-\sigma-1$ is at most
\begin{eqnarray}\label{eq:I}
&\ll&B^{|\sigma|}\max_{0\leq\nu\leq-\sigma-1}\frac{(|s|+\tau_k)^{2\nu}}{\nu!}\frac{(|s+\nu+i(\tau_k+2\theta_F)|+1)^{2|\sigma+\nu|+1}}{(|s|+\tau_k)^{2|\sigma|+1}}\nonumber \\
&\ll&B^{|\sigma|}\max_{0\leq\nu\leq-\sigma-1}\frac{(|\sigma|+\tau_k)^{2\nu}}{\nu!}
\frac{(|\sigma+\nu|+\tau_k)^{2|\sigma|-2\nu}}{(|\sigma|+\tau_k)^{2|\sigma|}}\\
&\ll& B^{|\sigma|}\max_{0\leq\nu\leq-\sigma-1}\frac{1}{\nu!}
\left(\frac{(|\sigma|-\nu+\tau_k)}{|\sigma|+\tau_k}\right)^{2(|\sigma|-\nu)} \ll B^{|\sigma|}.\nonumber
\end{eqnarray}
Still for $\si\leq -1$, the terms with $-\sigma\leq\nu\leq K$ (there are at most 3 of them) contribute
\begin{equation}\label{eq:II}
 \ll B^{|\sigma|}\frac{(|\sigma|+\tau_k)^{2|\sigma|}}{([|\sigma|]+1)!} \frac{\tau_k^A}{(|\sigma|+\tau_k)^{2|\sigma|}}\ll B^{|\sigma|}\frac{\tau_k^A}{(|\sigma|+1)^{|\sigma|}} \ll \tau_k^A.
 \end{equation}
 Finally, for $\si\leq -1$, the last term on the right hand side of \eqref{eq:0} contributes, recalling also \eqref{HK}, 
 \begin{equation}\label{eq:III}
 \frac{H_K(s+i\tau_k,p_a/q_F)}{F(s+i\tau_k)} \ll \
 B^{|\sigma|} \frac{|s+i\tau_k|^{2|\sigma|+A}}{(|\sigma|+\tau_k)^{2|\sigma|+1}}\ll 
 B^{|\sigma|} (|\sigma|+\tau_k)^A\ll B^{|\sigma|} \tau_k^A.
  \end{equation} 
 For $s\in\Omega_k$ and $\sigma\geq -1$ we have $|F(s+i\tau_k)|\gg \tau_k^{-A}$ whereas all the other terms in 
  (\ref{eq:0}) are $\ll \tau_k^A$. Therefore, gathering \eqref{eq:0},(\ref{eq:I}),(\ref{eq:II}) and (\ref{eq:III}) we conclude that 
 \begin{equation}\label{eq:IV}
 G_k(s)\ll B^{|\sigma|}\tau_k^A
 \quad \text{ for $s\in\Omega_k$.} 
 \end{equation}
 Note that the implied constant, as well as $A$ and $B$, may depend on $F$ but are independent of $\tau_k$.
 
 \smallskip
 Let now $s\in\Omega_{k,1}$. For such values we have $\tau_k^A\leq e^{A|\sigma|}$; hence, recalling our convention on the constant $B$, (\ref{eq:IV}) gives
 \begin{equation}\label{eq:V}
 G_k(s)\ll B^{|\sigma|} \qquad \text{for $s\in\Omega_{k,1}$.}
 \end{equation}
 For $s\in\Omega_{k,2}$ we write 
 \[
 g_k(s):=\exp(-4\cos(s/(2l_k))) \quad \text{and} \quad \widetilde{G}_k(s):=  B^s G_k(s)g_k(s). 
 \]
 We have $|g_k(s)|\leq 1$ for $s\in\Omega_{k,2}$, hence by (\ref{eq:Gk1}) we see that $\widetilde{G}_k(s)\ll 1$ for $s$ on the right vertical part of the boundary of $\Omega_{k,2}$. Similarly, using (\ref{eq:V}) we see that the same hold for $s$ on the left vertical part of the boundary of $\Omega_{k,2}$. On the horizontal parts of this boundary we have
$g_k(s)\ll\exp(-e^{l_k/2})$, therefore by (\ref{eq:IV}) we obtain that
\[
\widetilde{G}_k(s)\ll B^{\sigma+|\sigma|}\tau_k^A \exp(-e^{l_k/2})\ll 1.
\] 
Hence, thanks to the maximum modulus principle, we conclude that $\widetilde{G}_k(s)\ll 1$ for $s\in\Omega_{k,2}$. Suppose now that $s\in\Omega_{k,2}$ and $|t|\leq1$. Since $g_k(s)\gg\exp(-\exp(O(1/l_k))\gg 1$ and $\widetilde{G}_k(s)\ll 1$,
we have that $B^s G_k(s)\ll 1$ for such $s$. This, together with (\ref{eq:V}), implies that
\begin{equation}\label{eq:VI}
B^sG_k(s)\ll 1  \quad \text{for $\sigma\leq 1/2$ and $|t|\leq 1$}
\end{equation}
uniformly on $\tau_k$.

\smallskip
Now we recall that the functions $G_k(s)$ in \eqref{eq:Gkdef} were defined with the help of a sequence $\tau_k$ in Lemma \ref{lem:tauk}. In particular, we have $p^{-i\tau_k} \to \varepsilon_p$  for every $p|q_F$. Thus, in view of \eqref{effepi} and \eqref{eg:Gk}, for $\sigma>\vartheta'$ and $|t|\leq1$ the limit
\begin{equation}
\label{eq:IX}
G(s,{\boldsymbol \varepsilon}) = \lim_{k\to\infty} B^sG_k(s)
\end{equation}
exists and represents a bounded holomorphic function. Hence, by Vitali's convergence theorem (see Section 5.21 of Titchmarsh \cite{Tit/1952}) the limit exists and $G(s,{\boldsymbol \varepsilon})$ is holomorphic for all $\sigma\leq 1/2$, $|t|\leq 1$. Moreover, recalling (\ref{eq:VI}) and i) of Lemma 5 we have that
\begin{equation}\label{eq:X}
G(s,{\boldsymbol \varepsilon}) \ll 1
\end{equation}
uniformly  for ${\boldsymbol \varepsilon}\in {\mathbb T}^r$, $\sigma\leq 1/2$ and $|t|\leq 1$. Writing for $\sigma>\vartheta'$
\[
F_p(s)^{-1} = 1+ \sum_{m=1}^{\infty} \frac{c(p^m)}{p^{ms}},
\]
it is easy to see that
\[
H(s,\varepsilon_p):=\lim_{k\to\infty} \left(1-F_p(s+i\tau_k)^{-1}\right) = -\sum_{m=1}^{\infty} \frac{c(p^m)}{p^{ms}}\varepsilon_p^n.
\]
Therefore, recalling (\ref{eg:Gk}) and \eqref{eq:IX}, for $\sigma>\vartheta'$ we obtain
\begin{equation}\label{eq:XI}
G(s,{\boldsymbol \varepsilon})= B^s  \sum_{d|q_F}d\mu(\frac{q_F}{d}) \prod_{p|d} H(s,\varepsilon_p).
\end{equation}

\smallskip
Let now $d\mu(\boldsymbol \varepsilon)$ be the Haar measure on ${\mathbb T}^r=\prod_{p|q_F} {\mathbb T}_p$, where for every $p|q_F$ we denote by ${\mathbb T}_p$ a copy of the unit circle ${\mathbb T}$. Obviously, $d\mu(\boldsymbol \varepsilon)$ is the product of the normalized Lebesgue measures $d\mu(\ep_p)$ on the circles ${\mathbb T}_p$. For a given prime $p_j|q_F$ and a positive integer $m$ let
\[ 
J(s)=J(s;p_j,m):= \int_{{\mathbb T}^r} G(s,{\boldsymbol \varepsilon}) \varepsilon_{p_j}^{-m}\,
d\mu(\boldsymbol \varepsilon).
\]
Since $G(s,{\boldsymbol \varepsilon})$ is holomorphic for $\sigma\leq 1/2$ and $|t|\leq 1$, so is $J(s)$. Moreover, recalling  (\ref{eq:X}), for such $s$ we have
\begin{equation}\label{eq:J}
J(s)\ll  1.
\end{equation}
From (\ref{eq:XI}), Fubini's theorem and the orthogonality relations, for $\sigma>\vartheta'$  we have
\begin{eqnarray*}
J(s)&=& B^s  \sum_{d|q_F}d\mu(\frac{q_F}{d})  \int_{{\mathbb T}^r}\prod_{p|d} H(s,\varepsilon_p)  \varepsilon_{p_j}^{-m}\,
d\mu(\boldsymbol \varepsilon)\\
&=& B^s  \sum_{d|q_F}d\mu(\frac{q_F}{d}) \prod_{\substack{p|d\\ p\neq p_j}} \int_{{\mathbb T}_p}H(s,\varepsilon_p)\, 
d\mu(\varepsilon_p)  \int_{{\mathbb T}_{p_j}}H(s,\varepsilon_{p_j})  \varepsilon_{p_j}^{-m}\, 
d\mu(\varepsilon_{p_j})\\
&=& B^s p_j\mu(\frac{q_F}{p_j}) \int_{{\mathbb T}_{p_j}}H(s,\varepsilon_{p_j})  \varepsilon_{p_j}^{-m}\, 
d\mu(\varepsilon_{p_j}) = B^s \mu(\frac{q_F}{p_j}) c(p_j^m)p_j^{1-ms}.
\end{eqnarray*}
By analytic  continuation this equality holds for $\sigma\leq 1/2$ and $|t|\leq 1$, and using (\ref{eq:J}) we obtain
\[ 
c(p_j^m) \ll_{p_j,m} \left(\frac{B}{p_j^m}\right)^{|\sigma|}.
\]
For $m$ large enough the right hand side tends to $0$ as $\sigma\to-\infty$, thus $c(p_j^m)=0$ for such $m$. This shows that $F_{p_j}(s)^{-1}$ is a polynomial in $p_j^{-s}$, and Theorem 1 follows. \qed

\medskip
\section{Proof of Theorem 2} 

\smallskip
We follow closely the proof of Theorem 1.2 in \cite{Ka-Pe/2015}, to which we constantly refer. So here we shall be sketchy, indicating only the main changes. In the beginning we keep open the value of the sufficiently large constants $c_0,c_1, N$ below, and we add conditions when required. Moreover, we recall that the path $\LL(s)$, see Section 2, is defined as in Section 3.1 of  \cite{Ka-Pe/2015} but with a different choice of $t_0$; this implies some differences in the estimates below compared to the analogous estimates in Section 3.3 of \cite{Ka-Pe/2015}.

\smallskip
Let $z_X=\frac{1}{X}+2\pi i\alpha$ with a large $X>0$.  Writing
\[
F_X(s,\alpha) = \sum_{n=1}^\infty\frac{a(n)}{n^s}\exp(-nz_X),
\]
as in (3.27) of \cite{Ka-Pe/2015} for $\sigma<2$ we have
\begin{equation}
\label{P1}
F_X(s,\alpha) = \frac{1}{2\pi i}\int_{(2-\sigma)} F(s+w)\Gamma(w)z_X^{-w}\d w = \frac{1}{2\pi i}\int_{\LL(s)} F(s+w)\Gamma(w)z_X^{-w}\d w.
\end{equation}
If $w=u+iv\in\LL_{-\infty}(s)$ then $\Re(s+w)\geq -c_0$, hence 
$F(s+w)\ll |s+w|^c$ for some $c>0$ since $F(s)$ has polynomial growth on vertical strips. If in addition $v<t_0$, then $\Re(s+w)=c_0>1$ and hence $F(s+w)\ll 1$. Moreover, still for $w\in\LL_{-\infty}(s)$, we have
\[
|z_X^{-w}|=|z_X|^{-u}\exp(v(\pi/2-\eta(X))),
\]
where $\eta(X)>0$ and $\eta(X)=O(1/X)$, and by Stirling's formula
\[
\Gamma(w)  \ll e^{-\frac{\pi}{2}|v|}(|v|+1)^{u-\frac{1}{2}}.
\]
Therefore, due to the different choice of $t_0$, the contribution of $\LL_{-\infty}(s)$ to \eqref{P1} is
\[
\frac{1}{2\pi i}\int_{\LL_{-\infty}(s)} F(s+w)\Gamma(w)z_X^{-w}\d w  \ll A^{|\sigma|}(|s|+1)^{2|\sigma|+c}
\]
for some positive $A$ and $c$; such constants will not be necessarily the same at each occurrence. As a consequence, for $\sigma<2$ and any fixed $\alpha>0$, equation (3.28) of \cite{Ka-Pe/2015} becomes
\begin{equation}
\label{P2}
\begin{split}
F_X(s,\alpha) &= \frac{1}{2\pi i}\int_{\LL_{\infty}(s)} F(s+w)\Gamma(w)z_X^{-w}\d w + O(A^{|\sigma|}(|s|+1)^{2|\sigma|+c}) \\
&= \I_X(s,\alpha) + O(A^{|\sigma|}(|s|+1)^{2|\sigma|+c}),
\end{split}
\end{equation}
say, uniformly as $X\to\infty$. 

\smallskip
As on p.467 of \cite{Ka-Pe/2015}, now we apply the functional equation of $F(s)$ and the reflection formula of $\Gamma(s)$, thus getting
\[
\I_X(s,\alpha) = \omega Q^{1-2s} \frac{1}{2\pi i}\int_{\LL_{\infty}(s)} \overline{F}(1-s-w)G(s,w)S(s,w)(Q^2z_X)^{-w}\d w;
\]
see Section 2 for definitions. Taking into account the different choice of $t_0$, following Lemmas 3.1 and 3.2 of \cite{Ka-Pe/2015} for $w\in\LL_{\infty}$ we have
\[ 
S(s,w)=-ie(-\xi_F/4)e^{-\pi is} (1+ O(e^{-\eta v}))
\]
for a certain positive $\eta$ and
\[
G(s,w)\ll e^{-\frac{\pi}{2}(v+2t)}(v+|s|+1)^{2|\sigma|+c}.
\]
 Thus 
replacing $S(s,w)$ by $-ie(-\xi_F/4)e^{-\pi is}$ we obtain
\begin{equation}
\label{P3}
\begin{split}
\I_X(s,\alpha) &= -i\omega e(-\xi_F/4)Q^{1-2s}e^{-\pi is} \frac{1}{2\pi i}\int_{\LL_{\infty}(s)} \overline{F}(1-s-w)G(s,w)(Q^2z_X)^{-w}\d w \\
& \hskip1cm + O(A^{|\sigma|} \int_{t_0}^\infty|G(s,-\sigma-c_0+iv)| |(Q^2z_X)^{\sigma+c_0-iv}|
|e^{-\pi is}| e^{-\eta v}\d v) \\
&= \J_X(s,\alpha) + O(A^{|\sigma|}(|s|+1)^{2|\sigma|+c}),
\end{split}
\end{equation}
say, uniformly as $X\to\infty$.

\smallskip
As on p.467 of \cite{Ka-Pe/2015}, we reduce $G(s,w)$ in \eqref{P3} to a single $\Gamma$-factor by means of the uniform version of the Stirling formula in \cite{Ka-Pe/2011}. For $1 \leq N\leq |s|+c$ and $\beta= \prod_{j=1}^r\lambda_j^{2\lambda_j}$, arguing as in \cite{Ka-Pe/2015} we obtain
\begin{equation}
\label{P3bis}
\begin{split}
\log G(s,w) &= \log\Gamma(1-2s-w-2i\theta_F) + (\frac12-s-w)\log\beta + \log\prod_{j=1}^r\lambda_j^{-2i\Im\mu_j} \\
&+ \sum_{\nu=1}^N \frac{R_\nu(s)}{\nu(\nu+1)} \frac{1}{w^\nu} + O\big(\frac{(c(|s|+1))^{N+2}}{|w|^{N+1}}\big).
\end{split} 
\end{equation}
Since $O\big(\frac{(c(|s|+1))^{N+2}}{|w|^{N+1}}\big) < 1$ for $w\in\LL_\infty(s)$,  provided the constants in the definition of $\LL_\infty(s)$ are sufficiently large, for $1\leq N\leq |\sigma|+c$ we have 
\[
e^{O\big(\frac{(c(|s|+1))^{N+2}}{|w|^{N+1}}\big)} = 1 + O\big(A^{|\sigma|}\frac{(|s|+1)^{N+2}}{|w|^{N+1}}\big).
\]
Hence from \eqref{P3bis} we obtain
\begin{equation}
\label{P4}
\begin{split}
G(s,w) &= \Gamma(1-2s-w-2i\theta_F) \beta^{\frac12-s-w} \prod_{j=1}^r\lambda_j^{-2i\Im\mu_j} \exp\big(\sum_{\nu=1}^N \frac{R_\nu(s)}{\nu(\nu+1)} \frac{1}{w^\nu}\big) \\
&\times \big(1+ O\big(A^{|\sigma|}\frac{(|s|+1)^{N+2}}{|w|^{N+1}}\big)\big).
\end{split} 
\end{equation}
Moreover, using Lemma \ref{lemR}, for $w\in \LL_{\infty}(s)$ and $1\leq N\leq c(|s|+1)$ we have
\begin{equation}\label{**}
\begin{split}
\sum_{\nu=1}^N \frac{R_{\nu}(s)}{\nu(\nu+1)}\frac{1}{w^{\nu}} &\ll 
\sum_{\nu=1}^N \frac{(c'(|s|+1))^{\nu+1}}{\nu(\nu+1)}\frac{1}{|w|^{\nu}}\\
&\ll \sum_{\nu=1}^N \frac{1}{\nu(\nu+1)}\frac{(c'(|s|+1))^{\nu+1}}{(c_1(|s|+1))^{2\nu}}=
\sum_{\nu=1}^N \frac{1}{\nu(\nu+1)}(\frac{c'}{c_1})^{2\nu}\ll 1
\end{split}
\end{equation}
if $c_1\geq c'$. Again, we remark that \eqref{P4} and \eqref{**} hold thanks to our present choice of $t_0$ in the definition of $\LL_{\infty}(s)$. This small but significant change compared to \cite{Ka-Pe/2015} leads to much better estimates in the $s$-aspect (compare to Lemma 3.9 and (3.31) in \cite{Ka-Pe/2015}).

\smallskip
Next we replace $G(s,w)$ by its main term in \eqref{P4} inside the integral $\J_X(s,\alpha)$ in \eqref{P3}. This causes an error of size
\[
\begin{split}
&\ll A^{|\sigma|}e^{\pi t} \int_{\LL_\infty(s)} |\Gamma(1-2s-w-2i\theta_F)| \exp\big(\Re\big(\sum_{\nu=1}^N \frac{R_\nu(s)}{\nu(\nu+1)} \frac{1}{w^\nu}\big)\big) \big|(Q^2z_X)^{-w}\big| \frac{(|s|+1)^{N+2}}{|w|^{N+1}} |\d w|\\
&\ll A^{|\sigma|}e^{\pi t}(|s|+1)^{N+2} \int_{\LL_\infty(s)} |\Gamma(1-2s-w-2i\theta_F)| \big|(Q^2z_X)^{-w}\big|
\frac{|dw|}{|w|^{N+1}}\\
&\ll A^{|\sigma|}(|s|+1)^{N+2} \int_{t_0}^{\infty}e^{-\pi v/2} e^{\pi v/2} v^{-\sigma+c_0-N-1/2} \d v
\ll A^{|\sigma|}(|s|+1)^{N+2}
\end{split}
\]
if $N\geq -\sigma+c_0+1$, the bound being uniform in $X$. Indeed, the bound
\[
|\Gamma(1-2s-w-2i\theta_F) (Q^2z_X)^{-w}| |w|^{-N-1} \ll A^{|\si|} e^{-\pi t} v^{-\si + c_0-N-1/2},
\]
used to obtain the last estimate, follows by an application of Stirling's formula, namely 
\[
\begin{split}
\Gamma(1-2s-w-2i\theta_F) &\ll |1-2s-w-2i\theta_F|^{1/2-\si+c_0} \\
& \hskip.5cm \times \exp((2t+v+2\theta_F) \arg(1-2s-w-2i\theta_F)) \\
& \ll v^{1/2-\si+c_0} \exp((2t+v+2\theta_F) \arg(1-2s-w-2i\theta_F)),
\end{split}
\]
observing that
\[
\arg(1-2s-w-2i\theta_F) = -\arctan\left(\frac{2t+v+2\theta_F}{1-2\si-u}\right) =-\frac{\pi}{2} + O\left(\frac{|\si|+1}{2t+v+2\theta_F}\right).
\]
Therefore, from \eqref{P3}, \eqref{P4} and recalling the definition of $q_F$ and $\omega^*_F$, for 
\begin{equation}
\label{P5}
|\sigma|+c\leq N\leq |\sigma|+ c+1
 \end{equation}
 we have
\[
\begin{split}
\J_X(s,\alpha) &= \omega^*_F \big(\frac{q_F}{4\pi^2}\big)^{1/2 -s -i\theta_F} e(-\frac12(s+i\theta_F)) \frac{1}{2\pi i} \int_{\LL_{\infty}(s)} \overline{F}(1-s-w) \ \times \\
&\qquad \times \Gamma(1-2s-w-2i\theta_F) \exp\big(\sum_{\nu=1}^N \frac{R_\nu(s)}{\nu(\nu+1)} \frac{1}{w^\nu}\big)\big(\frac{q_Fz_X}{4\pi^2}\big)^{-w} \d w + O\big(A^{|\sigma|} (|s|+1)^{2|\sigma|+c'}\big).
\end{split}
\]
Hence, by the substitution $1-2s-w-2i\theta_F \to w$ in the above integral and recalling the definition of $\LL^*_{-\infty}(s)$ in Section 2 we obtain
\begin{equation}
\label{P6}
\begin{split}
\J_X&(s,\alpha) = -i\omega^*_F \big(\sqrt{q_F}\alpha - i\frac{\sqrt{q_F}}{2\pi X}\big)^{2s-1+2i\theta_F} \frac{1}{2\pi i} \int_{\LL^*_{-\infty}(s)} \bar{F}(s+w+2i\theta_F) \Gamma(w)  \\
& \times \exp\big(\sum_{\nu=1}^N \frac{(-1)^\nu R_\nu(s)}{\nu(\nu+1)} \frac{1}{(w+2s-1+2i\theta_F)^\nu}\big)\big(\frac{q_Fz_X}{4\pi^2}\big)^w \d w + O\big(A^{|\sigma|} (|s|+1)^{2|\sigma|+c'})\big),
\end{split}
\end{equation}
uniformly as $X\to\infty$.

\smallskip
 Writing again $\eta = 2s-1+2i\theta_F$, from the power series expansion of the exponential function and recalling \eqref{eqnew1} we have
\begin{equation}
\label{L3-15/1}
\begin{split}
\exp\big(\sum_{\nu=1}^N  \frac{(-1)^\nu R_\nu(s)}{\nu(\nu+1)} \frac{1}{(w+\eta)^\nu}\big) &= 1 + \sum_{\mu=1}^\infty \frac{V_{\mu,N}(s)}{(w+\eta)^\mu}\\
&=1 + \sum_{\mu=1}^N \frac{V_{\mu}(s)}{(w+\eta)^\mu} + O(\sum_{\mu=N+1}^{\infty} \frac{|V_{\mu,N}(s)|}{|w+\eta|^\mu}).
\end{split}
\end{equation}
Since for $w\in\LL^*_{-\infty}(s)$  we have $|w+\eta|\geq 2(c'(|s|+1))^2$, where $c'$ denotes the constant in Lemma \ref{lemVmuN}, using such a lemma we obtain
\[
\begin{split}
\sum_{\mu=N+1}^{\infty} \frac{|V_{\mu,N}(s)|}{|w+\eta|^\mu})&\ll
\sum_{\mu=N+1}^{\infty} \frac{(c'(|s|+1)^{2\mu}}{|w+\eta|^\mu}\leq 
 \frac{(c'(|s|+1)^{2(N+1)}}{|w+\eta|^{N+1}}\sum_{\mu=N+1}^{\infty} 2^{-\mu}\\
&\ll A^{|\sigma|} \frac{(|s|+1)^{2N+1}}{|w(w-1)\cdots(w-N)|}.
\end{split}
\]
Thus using \eqref{L3-15/1} and Lemmas \ref{lemVmuN} and \ref{lemwmu2}
\[
\begin{split}
&\exp\big(\sum_{\nu=1}^N  \frac{(-1)^\nu R_\nu(s)}{\nu(\nu+1)} \frac{1}{(w+\eta)^\nu}\big) =1 + \sum_{\mu=1}^N \frac{V_{\mu}(s)}{(w+\eta)^\mu}  + O(A^{|\sigma|} \frac{(|s|+1)^{2N+1}}{|w(w-1)\cdots(w-N)|})\\
&= 1+\sum_{\mu=1}^N V_\mu(s) \sum_{\nu=\mu}^N \frac{A_{\mu,\nu}(s)}{(w-1)\cdots(w-\nu)} 
+ O\big(\frac{4^NN!}{|w(w-1)\cdots(w-N)|} \sum_{\mu=1}^N |V_{\mu}(s)||\eta|^{N-\mu+1}\big)\\
 &\hspace{2.52in} + O(A^{|\sigma|} \frac{(|s|+1)^{2N+1}}{|w(w-1)\cdots(w-N)|})\\
 &= \sum_{\nu=0}^N \frac{Q_\nu(s)}{(w-1)\cdots(w-\nu)} + 
  O\big(\frac{4^NN!(c'(|s|+1)^N}{|w(w-1)\cdots(w-N)|} \sum_{\mu=1}^N \frac{(c'(|s|+1)^{\mu+1}}{(\mu-1)!}\big)\\
  &\hspace{1.8in}\ + O(A^{|\sigma|} \frac{(|s|+1)^{2N+1}}{|w(w-1)\cdots(w-N)|})\\
  &= \sum_{\nu=0}^N \frac{Q_\nu(s)}{(w-1)\cdots(w-\nu)} 
  + O(A^{|\sigma|} \frac{(|s|+1)^{2N+1}}{|w(w-1)\cdots(w-N)|})
\end{split}
\]
Now we replace the term $\exp(\sum_{\nu=1}^N...)$ in \eqref{P6} by the main term of the above formula. This  causes a further error of size
\[
\begin{split}
&\ll A^{|\sigma|}(|s|+1)^{2N+1} \int_{\LL^*_{-\infty}(s)} \frac{|\Gamma(w)|}{|w(w-1)\cdots(w-N)|} e^{\pi|v|/2} |\d w| \\
&= A^{|\sigma|}(|s|+1)^{2N+1} \int_{\LL^*_{-\infty}(s)} \frac{|\Gamma(w-N)|e^{\pi|v|/2}}{|w|} |\d w|.
\end{split}
\]
Moreover, we also want $N$ such that $\Re(w-N)\leq 0$ for $w\in\LL^*_{-\infty}(s)$, i.e. we choose
\begin{equation}
\label{P6bis}
N=[-\sigma]+k
\end{equation}
with a sufficiently large positive integer $k$ satisfying \eqref{P5}. With such a choice of $N$ we have $|\Gamma(w-N)|\ll e^{-\pi|v|/2}/|v|^{1/2}$, hence the integral is $\ll 1$. Therefore \eqref{P6} becomes, uniformly as $X\to\infty$,
\begin{equation}
\label{P7}
\begin{split}
\J_X(s,\alpha) &= -i\omega^*_F \big(\sqrt{q_F}\alpha - i\frac{\sqrt{q_F}}{2\pi X}\big)^{2s-1+2i\theta_F} \sum_{\nu=0}^N Q_\nu(s) \frac{1}{2\pi i} \int_{\LL^*_{-\infty}(s)} \overline{F}(s+w+2i\theta_F)  \ \times \\
&\qquad \times \Gamma(w-\nu) \big(\frac{q_Fz_X}{4\pi^2}\big)^w \d w + O\big(A^{|\sigma|} (|s|+1)^{2|\sigma|+B}),
\end{split}
\end{equation}
which is the analog of (3.35) in \cite{Ka-Pe/2015}.

\smallskip
Replacing the path of integration in \eqref{P7} by the whole path $\LL^*(s)$ causes an error which, since $\nu\leq N$, by Lemma \ref{3.18} is of size
\[
\ll A^{|\sigma|} \sum_{\nu=0}^N \frac{(|s|+1)^{2\nu}}{\nu!} \int_{\LL_\infty^*(s)} |F(s+w+2i\theta_F) \Gamma(w-\nu) \big(\frac{q_Fz_X}{4\pi^2}\big)^w| |\d w|.
\]
Let now $w\in\LL^*_\infty(s)$. Clearly, $F(s+w+2i\theta_F)\ll 1$. Moreover,  for $0\leq \nu\leq N$ we have $-c \leq \Re(w-\nu) \leq N+1-\nu$ and hence for $v\neq0$
\[
v\arg(w-\nu) = v \arg(iv) + \arg\left(1+i\frac{\nu-u}{v}\right) = \frac{\pi}{2}|v| 
+O(|\si|+1).
\]
Observing that this estimate holds for $v=0$ as well, by Stirling's formula we get
\[
\log|\Gamma(w-\nu)| = (u-\nu-1/2)\log|w-\nu| -\frac{\pi}{2}|v| + O(|\si|+1).
\]
Consequently
\[
\Gamma(w-\nu) \ll A^{|\si|} e^{-\frac{\pi}{2}|v|} (|s|+|v|+1)^{N-\nu+1/2}.
\]
Further, $|z_X^w| \ll A^{|\sigma|} e^{-v\arg z_X}$, thus the above mentioned error is (with a suitable $c''>0$)
\[
\ll A^{|\sigma|} \sum_{\nu=0}^N \frac{(|s|+1)^{2\nu}}{\nu!} \int_{-t_0^*(s)}^\infty (|s|+|v|+1)^{|\si|-\nu+c''}
e^{-\pi|v|/2-v\arg{z_X}}\d v \ll A^{|\sigma|}(|s|+1)^{2|\sigma|+B}.
\]
 We also have by Cauchy's theorem that for $0\leq \nu\leq N$
\[
\begin{split}
&\frac{1}{2\pi i} \int_{\LL^*(s)} \overline{F}(s+w+2i\theta_F) \Gamma(w-\nu) \big(\frac{q_Fz_X}{4\pi^2}\big)^w \d w \\
&\hskip-.5cm = \frac{1}{2\pi i} \int_{|\sigma|+\nu+2-i\infty}^{|\sigma|+\nu+2+i\infty}  \overline{F}(s+w+2i\theta_F) \Gamma(w-\nu) \big(\frac{q_Fz_X}{4\pi^2}\big)^w \d w \\
&\hskip-.5cm =  \big(\frac{q_F}{4\pi^2X}+i\frac{q_F\alpha}{2\pi}\big)^\nu \sum_{n=1}^\infty \frac{\overline{a(n)}}{n^{s+\nu+2i\theta_F}} \exp(-\frac{4\pi^2}{q_Fz_X}n)
\end{split}
\]
since the poles of the integrand lie to the left of $\LL^*(s)$.
Consequently, \eqref{P7} becomes
\begin{equation}
\label{P8}
\begin{split}
\J_X(s,\alpha) &= -i\omega^*_F \big(\sqrt{q_F}\alpha - i\frac{\sqrt{q_F}}{2\pi X}\big)^{2s-1+2i\theta_F} \sum_{\nu=0}^N 
 \big(\frac{q_F}{4\pi^2X}+i\frac{q_F\alpha}{2\pi}\big)^\nu Q_\nu(s) \\
&\hskip.8cm \times \sum_{n=1}^\infty \frac{\overline{a(n)}}{n^{s+\nu+2i\theta_F}} \exp(-\frac{4\pi^2}{q_Fz_X}n) +O(A^{|\sigma|}(|s|+1)^{2|\sigma|+B}),
\end{split}
\end{equation}
uniformly as $X\to\infty$. Note that, since
\[
-\frac{4\pi^2}{q_Fz_X} = \frac{2\pi i}{q_F\alpha} - \frac{1}{q_F\alpha^2}\frac{1}{X+O(1)},
\]
the series in \eqref{P8} is absolutely convergent for all $s$, for every $\nu$.

\smallskip
As on p.472 of \cite{Ka-Pe/2015}, the final step is to make the range of summation of $\nu$ in \eqref{P8} independent of $\sigma$ (recall that $N$ depends on $\sigma$, see \eqref{P6bis}). Let $K>0$ be a large integer and  $\sigma> -K+1/2$. Depending on the relative sizes of $N$ and $K$, we add to or withdraw from \eqref{P8} the terms with $\nu$ between $N+1$ and $K$ or between $K+1$ and $N$, respectively. 
In both cases we have that $\sigma+\nu> 3/2$ for such $\nu$'s (call them $\nu\in\X$), hence from Lemma 
\ref{3.18} we deduce that
\begin{equation}
\label{P9}
\begin{split}
-i\omega^*_F &\big(\sqrt{q_F}\alpha - i\frac{\sqrt{q_F}}{2\pi X}\big)^{2s-1+2i\theta_F} \sum_{\nu\in\X} \big(\frac{q_F}{4\pi^2X}+i\frac{q_F\alpha}{2\pi}\big)^\nu Q_\nu(s) \times \\
&\times \sum_{n=1}^\infty \frac{\overline{a(n)}}{n^{s+\nu+2i\theta_F}} \exp(-\frac{4\pi^2}{q_Fz_X}n) 
\ll A^{|\sigma|} \sum_{\nu\in\X} \frac{(c'(|s|+1))^{2\nu}}{\nu!} 
\end{split}
\end{equation}
uniformly in $X$. If $N<K$ this is
\[
\ll A^{|\sigma|}(|s|+1)^{2K} \sum_{\nu\geq N+1} \frac{c'^{2\nu}}{\nu!} \ll
\frac{ A^{|\sigma|}c'^{2N}}{(N+1)!}(|s|+1)^{2K} \sum_{\nu\geq0} \frac{c'^{2\nu}}{\nu!} 
\ll (|s|+1)^{2K},
\]
while if $N>K$ we have 
\[
A^{|\sigma|} \sum_{\nu\in\X} \frac{(c'(|s|+1))^{2\nu}}{\nu!} \ll 
A^{|\sigma|} \frac{(c'(|s|+1)^{2N}}{(N-1)!} \ll (|s|+1)^{2N}\ll (|s|+1)^{2K+A}.
\]
Thus, in view of \eqref{P2},\eqref{P3},\eqref{P8} and \eqref{P9}, for $-K+1/2<\sigma<2$ we have
\begin{equation}
\label{P10}
\begin{split}
F_X(s,\alpha) = &-i\omega^*_F \big(\sqrt{q_F}\alpha - i\frac{\sqrt{q_F}}{2\pi X}\big)^{2s-1+2i\theta_F} \sum_{\nu=0}^K \big(\frac{q_F}{4\pi^2X}+i\frac{q_F\alpha}{2\pi}\big)^\nu \\
&\times Q_\nu(s) F^*_X(s+\nu+2i\theta_F,\alpha) + H_X(s,\alpha),
\end{split}
\end{equation}
where
\[
F^*_X(s,\alpha) = \sum_{n=1}^\infty \frac{\overline{a(n)}}{n^{s}} \exp(-\frac{4\pi^2}{q_Fz_X}n)
\]
and 
\[
H_X(s,\alpha) \ll (|s|+1)^{2K+A}
\]
uniformly as $X\to\infty$. Moreover, since $F_X(s,\alpha)$, $F^*_X(s,\alpha)$ and $Q_\nu(s)$ are entire functions, $H_X(s,\alpha)$ is also entire. Further, from \eqref{P10} we have that for $1<\sigma<2$
\[
\lim_{X\to\infty} H_X(s,\alpha) = H(s,\alpha)
\]
exists and is holomorphic since this is clearly true for $F_X(s,\alpha)$ and $F^*_X(s,\alpha)$. Thanks to \eqref{P10}, for $1<\sigma<2$
we also have that 
\begin{equation}
\label{P11}
H(s,\alpha) = F(s,\alpha) + i\omega^*_F \big(\sqrt{q_F}\alpha\big)^{2s-1+2i\theta_F} \sum_{\nu=0}^K \big(i\frac{q_F\alpha}{2\pi}\big)^\nu Q_\nu(s) \bar{F}(s+\nu+2i\theta_F,-\frac{1}{q_F\alpha}).
\end{equation}
Hence by Vitali's convergence theorem the limit function $H(s,\alpha)$ exists and is holomorphic for $-K+1/2<\sigma<2$, and satisfies
\[
H(s,\alpha) \ll (|s|+1)^{2K+A}.
\]
This provides analytic continuation and bounds for the right-hand side of \eqref{P11}, and Theorem 2 follows.

\bigskip
\bigskip

\bigskip
\bigskip
\noindent
Jerzy Kaczorowski, Faculty of Mathematics and Computer Science, A.Mickiewicz University, 61-614 Pozna\'n, Poland and Institute of Mathematics of the Polish Academy of Sciences, 00-956 Warsaw, Poland. e-mail: \url{kjerzy@amu.edu.pl}

\medskip
\noindent
Alberto Perelli, Dipartimento di Matematica, Universit\`a di Genova, via Dodecaneso 35, 16146 Genova, Italy. e-mail: \url{perelli@dima.unige.it}

\end{document}